\documentclass{article}
\usepackage{listings}
\usepackage{xcolor}
\usepackage{verbatim}
\usepackage{latexsym}
\usepackage{dsfont}
\usepackage{amssymb}
\usepackage{tocbibind}
\usepackage{amsthm}
\usepackage{amsmath}
\usepackage{mathtools}
\usepackage{makecell}
\usepackage{graphicx}
\usepackage{tikz}
\usepackage{mathrsfs}
\usepackage{cite}
\usepackage{appendix}
\usepackage{indentfirst}

\def\be{\begin{eqnarray}}
\def\ee{\end{eqnarray}}
\def\X{{\cal X}}

\def\({\left(}
\def\){\right)}
\def\bc{\begin{center}}

\def\ec{\end{center}}

\def\bey{\begin{eqnarray*}}
\def\eey{\end{eqnarray*}}
\def\ber{\begin{array}{l}}
\def\eer{\end{array}}
\newtheorem{theorem}{Theorem}

\newtheorem{corollary}{Corollary}
\newtheorem{definition}{Definition}

\begin{document}
\title{Some connections between classical and nonclassical symmetries of a partial differential equation and their applications\footnote{Supposed by NSF of China(No.11571008)}}


\author{Temuer Chaolu$^1$, Laga Tong$^2$, George Bluman$^3$\\
{\small $^1$Shanghai Maritime University~~Shanghai~~China~~201306}\\
{\small $^2$\footnote{$^2$Corresponding author, e-mail: lagelage@163.com}Inner Mongolia University of Technology~~Huhhot~~China~~010051}\\
{\small $^3$University of British Columbia~~Vancouver~~BC~~
Canada V6T 1Z2}
\\}

\date{}
\maketitle

{\bf Abstract: } Some connections between classical and nonclassical symmetries of a partial differential equation (PDE) are given in terms of determining equations of the two symmetries. These connections provide additional information for determining nonclassical symmetry of a PDE and make it easier to solve the system of nonlinear determining equations. As example, new nonclassical symmetries are exhibited for a class of generalized Burgers equation and KdV-type equations are given.


\section{Introduction}
The nonclassical symmetry method, proposed by Bluman \textit{et al} in \cite{B2}, is one of generalization of the classical Lie symmetry  method for obtaining solutions of nonlinear PDEs \cite{b3, o1}. The nonclassical symmetry admitted by a PDE yields exact solutions to the PDE which can not be obtained through its classical symmetry. It is also related to B{$\rm\ddot{a}$}cklund transformations, functionally invariant solutions, and the ¡±direct¡± method etc \cite{G, Nucci, Gerd}.

Determining classical and nonclassical symmetries of a PDE is equivalent to exactly solving the so-called the systems of determining equations arising from the invariance of the PDE under such symmetry transformations. There exist algorithms and software packages  \cite{Man2, D,C, chao5, chao2} to  produce such systems of determining equations. However, to obtain nonclassical symmetries, exactly solving the system is a hard task due to its nonlinearity unlike the situation for classical symmetry \cite{A4}. For this reason, nonclassical symmetries of many PDEs in physics and mechanics have not been found \cite{I}. In \cite{D, C, chao5, chao3, chao2}, the authors investigate ways of simplifying the solution of the system of nonlinear determining equations based on the Gr${\rm\ddot{o}}$ebner basis method and Wu's method. In the both methods, a system of determining equations is regarded as one of differential polynomials. Here, the set of solutions (equivalent to the set of symmetries of the PDE) of the system of determining equations is the same as the set of zero points (algebra variety) of the corresponding system of differential polynomials. Consequently, the problem of solving the system of determining equations is turned to determining the zero set of the corresponding polynomial system. This is suitably dealt with within the frameworks of  the Gr${\rm\ddot{o}}$ebner basis method and Wu's method. These methods lead to the development of algorithms for directly dealing with the system of determining equations. However, in these algorithms the computation of the differential Gr${\rm \ddot{o}}$ebner basis or the characteristic set of the system of corresponding determining equations is a difficult task and this is hardly implemented directly for higher dimensional PDEs\cite{A4}.

The aim of this article is to further investigate the connections between the classical and nonclassical symmetries by establishing relationships between the systems of determining equations for these symmetries of a PDE. Particularly, we will show the connections through comparing the determining polynomials of the nonclassical symmetry with those of the classical symmetry. This leads to the essential links between classical and nonclassical symmetries of a PDE in terms of their respective determining equations.

An innovation in this article is that some information on the nonclassical symmetries of a PDE can be observed directly from its classical symmetries. In particular, we provide special supplementary conditions for nonclassical symmetries so that the nonlinear determining system is easier to solve. Our results not only give an alternative way to determine nonclassical symmetries of a PDE, but also improve the efficiency of the method in \cite{chao5}.

The rest of this paper is organized as follows. In Section 2, we give some preliminaries on the classical and nonclassical symmetries of a PDE and present basic results of Ritt-Wu's method used in this article. In Section 3, we present our main results on connections between classical and nonclassical symmetries of a PDE in terms of the determining polynomial systems of the symmetries. In Section 4, as applications of the obtained connections, new nontrivial nonclassical symmetries are found for both a class of generalized Burgers equations and system of KdV-type equations. In Section 5, we summarize our results and include concluding remarks.

\section{Preliminary}
In the following, we give some preliminaries used in this article.
\subsection{Notations}

Let $\mathbf{x}=(x_1, x_2, \ldots, x_p)$ and $\zeta=\{\xi_1, \xi_2, \ldots, \xi_q\}$ be independent and dependent variables, respectively. In terms of a multi-index $\alpha=\{\alpha_1, \ldots, \alpha_p \} \in \mathbb{N}_0^p$($\mathbb{N}_0$ is the set of non-negative integers), denote the total differential operator $D^\alpha=\partial^{|\alpha|}/\partial_{x_1}^{\alpha_1}\partial_{x_2}^{\alpha_2}\ldots$ $ \partial_{x_p}^{\alpha_p}$ with order $|\alpha|=\sum_{i=1}^p\alpha_i$; for an $k\in \mathbb{N}_0$, we use the notation $\partial^k \zeta=\{D^\alpha\xi_i, \alpha\in \mathbb{N}_0^p,|\alpha|\leq k, i=1, 2, \ldots, q\}$. Let $\mathcal{K}_\mathbf{x}$ be a characteristic zero field of differential functions of $\mathbf{x}$ and let $\mathcal{K}_\mathbf{x}[\zeta]$ be the differential polynomial ring in indeterminates $\partial^k\zeta$ for all $k\in \mathbb{N}_0$ over $\mathcal{K}_\mathbf{x}$.

A system of differential equations $\{p_1=0, p_2=0,\ldots, p_s=0\}$ where $p_i$ are polynomials in $\mathcal{K}_\mathbf{x}[\zeta]$ is denoted shortly as $\mathcal{D}=0$ where $\mathcal{D}$ is the set of left side polynomials of the equations, i.e., $\mathcal{D}=\{p_1, p_2, \ldots, p_s\}\subset \mathcal{K}_\mathbf{x}[\zeta]$. For a system of differential equations $\mathcal{D}=0$ in polynomial form, one has a corresponding system of differential polynomials $\mathcal{D}$.
\rm{Z}$(\mathcal{D})$  denotes the set of zero points of a differential polynomial system $\mathcal{D}$ over some extended field of $\mathcal{K}_\mathbf{x}$. Thus the solution set to the system of equations $\mathcal{D}=0$ is the same as the set Z$(\mathcal{D})$.

For two systems $\mathcal{D}$ and $Q$ of differential polynomials, we denote $\textrm{Z}(\mathcal{D}/Q)=\textrm{Z}(\mathcal{D})\backslash\textrm{Z}(Q)$.

The notation $f|_{\textrm{Z}(\mathcal{D})}$ denotes the set of values of the  differential polynomial $f\in \mathcal{K}_\mathbf{x}[\zeta]$ at all points in Z$(\mathcal{D})$.

In all examples, we use the graded lexicographic order \cite{chao5, D, ritt, kolch} for differential polynomials with no further elaboration.

\subsection{A reduction formula}

A reduction of a differential polynomial with respect to another one or an ascending chain of differential polynomials is a basic operation in polynomial algebra. The frequently used reduction formula is the pseudo-reduction given in \cite{ritt, kolch}. It has been used in the methods given in \cite{D, Man2} and \cite{wu3, wu4} under different names such as Kolchin-Ritt reduction and Ritt-Wu reduction. We apply the reduction in this article for setting up a connection between classical and nonclassical symmetries of a PDE and analyze the set of zero points of a system of polynomials.


The concept of a differential ascending chain of a differential polynomial system is given below.
\begin{definition}A finite differential polynomial system:
\begin{eqnarray} \mathcal{AS}=\{\mathcal{A}_1, \mathcal{A}_2, \ldots, \mathcal{A}_s\}\subset \mathcal{K}_\mathbf{x}[\zeta], \label{dchain}
\end{eqnarray}
is called a differential ascending chain under a differential polynomial rank $\prec$ if it satisfies the following two conditions:

{\rm(a)} $\mathcal{A}_1\prec \mathcal{A}_2\prec \ldots \prec\mathcal{A}_s$;

{\rm(b)} $\mathcal{A}_j$ is reduced with respect to $\mathcal{A}_i$ {\rm(}the terms in $\mathcal{A}_j$ can not be eliminated by $\mathcal{A}_i${\rm )} for $i=1, 2, \ldots, j-1.$
\end{definition}

The Ritt-Wu reduction for a differential polynomial $f\in \mathcal{K}_\mathbf{x}[\zeta]$ with respect to a differential ascending chain $\mathcal{AS}\subset\mathcal{K}_\mathbf{x}[\zeta]$ is given as follows.

\begin{theorem}\label{th0}For a differential polynomial $f\in \mathcal{K}_\mathbf{x}[\zeta]$ and a differential ascending chain $\mathcal{AS}\subset \mathcal{K}_\mathbf{x}[\zeta]$, there exist differential operators $\mathfrak{D}_i=\sum_{\beta\in \mathds{N}_0^p}{Q_i^\beta}D_i^\beta$ with $Q_i^\beta\in \mathcal{K}_\mathbf{x}[\zeta]$ such that
\begin{eqnarray} {\rm IS}\cdot f=\sum_{g_i\in
\mathcal{AS}} \mathfrak{D}_i g_i+r,\label{redf}
\end{eqnarray}
where the differential polynomial $r\in \mathcal{K}_\mathbf{x}[\zeta]$ is reduced with respect to $\mathcal{AS}$ {\rm (}$r$ cannot be reduced further by any polynomial in $\mathcal{AS}$ and is called the pseudo-remainder of $f$ with respect to $\mathcal{AS}$, denoted by $r ={\rm Prem(}f/\mathcal{AS}{\rm )}$. The $D_i^\beta$ is the total differential operator related to $g_i\in \mathcal{AS}$ with respect to $\beta\in \mathds{N}^p_0$ introduced in Section 2.1. The {\rm IS (}called IS product of $\mathcal{AS}${\rm )} is a product of initials {\rm (}the coefficients of the leading derivatives of polynomials{\rm )} and separants {\rm (}partial derivatives of differential polynomials with respect to leading derivatives{\rm )}  of the differential polynomials in $\mathcal{AS}$.
\end{theorem}
The reduction formula (\ref{redf}) is proved by successively eliminating from $f$ all derivative terms that can be obtained by differentiation of the leading derivatives of the  elements of $\mathcal{AS}$ under a polynomial order until the elimination cannot be continued (the details are seen in \cite{wu3}). The procedure of building (\ref{redf}) is constructive and can be implemented by a computer algebra software system \cite{chao5}. For a differential polynomial system $\mathcal{D}$, we use the notation
\begin{eqnarray*}
{\rm Prem(}\mathcal{D}/\mathcal{AS}{\rm ) =} \{{\rm Prem(}f/\mathcal{AS})\textrm{  for } f \in
\mathcal{D}\}.
\end{eqnarray*}

In addition, a basic algorithm used in this article is the well-known Wu's algorithm (also called the characteristic set algorithm) for constructing a differential characteristic set  for any given finite differential polynomial system. We use the algorithm to reduce a system of determining equations of the symmetries of a PDE into a differential ascending chain (see the given algorithm and examples in the following sections). The details of the algorithm appear in \cite{wu3, wu4}.

\subsection{Classical and nonclassical symmetries of a PDE}
To make our statement more concise, we consider the case of a scale PDE. Moreover, the results presented in the article hold in the case of a system of PDEs.

The general statement on the nonclassical symmetry for a general PDEs appears in \cite{D}.

Let an $k$th order PDE
\begin{eqnarray}
{\rm F}[u]={\rm F}\left(\mathbf{x},  \partial^k u\right)=0, \label{e1}
\end{eqnarray}
be given with independent variables $\mathbf{x}=(x_1, x_2, \ldots, x_p)$ and one dependent variable $u$. Suppose $\textrm{F}$ is a polynomial in its arguments.

Let
\begin{eqnarray}
\X^\prime=\sum_{i=1}^p\xi_i^\prime\partial_{x_i}+\eta^\prime\partial_{u}, \label{e3}
\end{eqnarray}
be the infinitesimal generator of classical symmetry of equation (\ref{e1}) where $\xi_i^\prime=\xi_i^\prime(\mathbf{z})$ and $\eta^\prime=\eta^\prime(\mathbf{z})$ are the infinitesimal functions of the generator with independent variables $\mathbf{z}=(x_1, x_2, \ldots, x_p, u)$.

In the classical symmetry method, the system of determining equations of generator (\ref{e3}), a system of over-determined homogeneous linear PDEs for the functions $\xi_i^\prime(\mathbf{z})$ and $\eta^\prime(\mathbf{z})$, is derived from the invariance criterion  ${\rm{Pr}}\X^\prime({\rm F})=0$ when ${\rm F}=0$. Here $\rm{Pr}\mathcal{X}^\prime$ is the prolongation of $\mathcal{X}^\prime$ on the jet space of $u$. Explicitly solving the system of determining equations, one determines the generator (\ref{e3}), i.e., the classical symmetry of PDE (\ref{e1}).

In the standard nonclassical symmetry method\cite{D,B2}, the infinitesimal generator
\begin{equation}
\mathcal{X}=\sum_{i=1}^p\xi_i\partial_{x_i}+\eta\partial_{u},\label{ninf}
\end{equation}
of the nonclassical symmetry of equation (\ref{e1}) is found by requiring that both equation (\ref{e1}) and its invariance surface equation
\begin{eqnarray}
\psi =\sum_{i=1}^p\xi_i u_{x_i}-\eta =0,\label{e4}
\end{eqnarray}
are simultaneously invariant under the symmetry transformation of (\ref{ninf}). In this case, the invariance criterion is given by $\rm{Pr}\X(F)=0$ which holds on the manifold defined by (\ref{e1}) and the differential consequences of (\ref{e4}). From this criterion, one derives the system of determining equations for a nonclassical symmetry (\ref{ninf}) which is a nonlinear system of PDEs for the infinitesimal functions $\xi_i=\xi_i(\mathbf{z})$ and $\eta=\eta(\mathbf{z})$ that appear in (\ref{ninf}) and  (\ref{e4}). Thus, in principle, after exactly solving the system, one finds the infinitesimal generator (\ref{ninf}), i.e., the nonclassical symmetry of system (\ref{e1}). However, this is, in general, a difficult task \cite{D, A4}.

The concept of equivalent generators plays a key rule in this article.

\begin{definition}Two first order differential operators $\mathcal{X}^\prime$ and $\mathcal{X}$ of the form as {\rm(\ref{e3})} and {\rm (\ref{ninf}) }are called equivalent if they differ by a nonzero differential function multiplier $\lambda(\textbf{z})$, i.e., $\mathcal{X}^\prime=\lambda(\textbf{z})\mathcal{X}$.
\end{definition}
Equivalent symmetry generators of a PDE generate the same reduced solutions of the PDE.

\textbf{Remark 1:} In this article, without loss of generality, we assume $\xi_1\not = 0$ in (\ref{ninf}). Otherwise, we successively assume $\xi_k\not=0$ and $\xi_i=0, i=1, 2, \cdots, k-1$ and discuss each case separately. Hence, by Definition 2, we set $\xi_1=1$. Conventionally, in the case of an evolution equation case, each $\xi_1$ is the coefficient of $\partial_t$ in generator $\mathcal{X}$ corresponding to the time variable, usually denoted by $\tau$. Hence, in the regular case of the nonclassical symmetries for an evaluation equation, we set $\tau=1$.

In the following, we let $\mathcal{D}^\prime=0$ and $\mathcal{D}=0$ be the systems of determining equations for the symmetry generators $\mathcal{X}^\prime$ and $\mathcal{X}$ of PDE (\ref{e1}). The sets of the left hand sides $\mathcal{D}^\prime\subset \mathcal{K}_{\mathbf{z}}[\Lambda^\prime]$ and $\mathcal{D}\subset \mathcal{K}_{\mathbf{z}}[\Lambda]$ of the two systems of determining equations are called determining polynomial systems of the classical and nonclassical symmetries respectively. Consequently, the zero sets ${\rm{Z}}(\mathcal{D}^\prime)$ and ${\rm{Z}}(\mathcal{D})$ represent the solution sets of the two systems of determining equations, i.e., the sets of the classical and nonclassical symmetries of the PDE respectively.

For the sake of simplicity, we use following additional notations: $\Lambda^\prime=(\xi_1^\prime, \ldots, \xi_p^\prime,$ $\eta^\prime)$,  $\Lambda=(\xi_2, \ldots, \xi_p, \eta)$, $\partial_\textbf{z}=(\partial_{x_1}, \partial_{x_2}, \ldots, \partial_{x_p}, \partial_{u})$ and $\tilde{\partial}_\textbf{z}=(\partial_{x_2}, $ $ \partial_{x_3}, \ldots, \partial_{x_p},$ $\partial_{u})$.
Thus, the generators (\ref{e3}) and (\ref{ninf}) can be written shortly as $\mathcal{X}^\prime=\Lambda^\prime\cdot \partial_\textbf{z}$ and $\mathcal{X}=\partial_{x_1}+\Lambda\cdot \tilde{\partial}_\textbf{z}$ by scalar product rule of finite dimensional vector space.

Our main goal is to determine some connections between the systems $\mathcal{D}$ and $\mathcal{D}^\prime$.

\section{Some connections between classical and nonclassical symmetries of a PDE}

It is well-known that the set of classical symmetries of a PDE is a subset of the nonclassical symmetries of the PDE. Hence, the classical symmetry of a PDE has equivalent expressions in the coordinates $\Lambda^\prime$ and $\Lambda$.

Therefore, we use a map between any $\Lambda^\prime=(\xi_1^\prime, \cdots, \xi_p^\prime, \eta^\prime)$ with $\xi_1^\prime\not=0$ and $\Lambda=(\xi_2, \cdots, \xi_p, \eta)$ as coordinately
\begin{eqnarray}
\xi^\prime_i=\xi_1^\prime \xi_i,\,\,\,\eta^\prime=\xi_1^\prime\eta,\,\,\, i=2, \ldots, p. \label{CLN1}
\end{eqnarray}

If a nonclassical symmetry $\mathcal{X}=\partial_{x_1}+\Lambda \cdot \tilde{\partial}_\textbf{z}$ represents a classical symmetry, then there exists an equivalent $\mathcal{X}^\prime=\Lambda^\prime\cdot\partial_\textbf{z}$ such that $\Lambda^\prime=(\xi_1^\prime, \cdots, \xi_p^\prime, \eta^\prime)\in \textrm{Z}(\mathcal{D}^\prime)$  and $(1,\Lambda)=(1,\xi_2, \cdots, \xi_p, \eta)\in \textrm{Z}(\mathcal{D})$ have a unified proportional relationship. Hence, the transformations (\ref{CLN1}) map a classical symmetry expressed by an $\Lambda\in {\rm Z}(\mathcal{D})$ into equivalence one expressed by the corresponding $\Lambda^\prime \in {\rm Z}(\mathcal{D}^\prime)$ and vice versa.

We take a differential polynomial rank $\prec$ on  $\mathcal{K}_{\mathbf{z}}[\Lambda^\prime]$ with $\xi^\prime_1$ having the highest rank, i.e., $\xi_j^\prime\prec \xi_1^\prime, j\not=1$.  Let $C^\prime$ be an ascending chain form of $\mathcal{D}^\prime$ which can be obtained by Wu's algorithm under the rank. In  \cite{chao5}, we have proved that ${\rm{Z}}(C^\prime)={\rm{Z}}(\mathcal{D}^\prime)$ for a PDE if ${\rm IS}(C^\prime)\not=0$. Hence, the set ${\rm{Z}}(C^\prime)$ represents the classical symmetry of the PDE under the conditions that the initials of $C^\prime$ are not zero.

\textbf{Remark 2: }Here, we emphasize that it is unnecessary for $C^\prime$  to be a characteristic set of $\mathcal{D}^\prime$, which is required in the algorithm given in \cite{chao5}. The cost of calculating $C^\prime$ is much less than that of calculating the characteristic set of $D^\prime$. Particularly, it is noted that the set $C^\prime$ is an irreducible ascending chain due to its linearity.

In this article, we assume that the initials of the differential polynomials in $C^\prime$ are not zero, i.e.,
\begin{eqnarray}
{\rm{IS}(}C^\prime{\rm)}\not\equiv 0.\label{init}
\end{eqnarray}

\textbf{Remark 3: }Assumption (\ref{init}) is not critical for a PDE with no arbitrary parameters. Since, in this case, the initial product IS of $C^\prime$ is a polynomial in terms of the independent variable $\textbf{z}$, so that it is never zero.

Now, we can obtain the first connection between classical and nonclassical symmetries of a PDE, which allows one to determine the existence of a nontrivial nonclassical symmetry.
\begin{theorem}\label{th2}
For an $\Lambda \in{\rm Z}(\mathcal{D})$, if, after using transformations {\rm (\ref{CLN1})} for a $\xi_1^\prime\not=0$, the $C^\prime$ as a differential polynomial system in $\mathcal{K}_{\bf z}[\Lambda, \xi_1^\prime]$ has no zero points, i.e., {\rm Z}$(\mathcal{C}^{\prime})=\emptyset$, then the operator $\mathcal{X}=\partial_{x_1}+\Lambda\cdot \tilde{\partial}_{\mathbf{z}}$ is a nontrivial nonclassical symmetry of the \emph{PDE}.
\end{theorem}
Proof: If the $\mathcal{X}$ is equivalent to a classical symmetry, then there exists  $\Lambda^\prime=(\xi_1^\prime,\cdots, \xi_p^\prime, \eta^\prime)\in {\rm Z}(C^\prime)$ under transformation (\ref{CLN1}). This means $\textrm{Z}(C^\prime)\not=\emptyset$. This is a contradiction with respect to the assumption of the theorem. The theorem is proved.

In particular, suppose $\mathcal{D}^{\prime\prime}$ be the subset of $C^\prime$ consisting of the differential polynomials whose leading derivatives are the derivatives of $\xi_1^\prime$ with respect to $\mathbf{z}$. For any subset $Q\subset C^\prime$, one has the following from this theorem.

\begin{corollary}
The theorem \ref{th2} is also true when replacing $C^\prime$ in the theorem with any subset $Q\subset C^\prime$. In particular, when $Q$ is $\mathcal{D}^{\prime\prime}$ or the complement set $C^\prime\backslash \mathcal{D}^{\prime\prime}$, the theorem is true.
\end{corollary}

In the following, based on formula (\ref{redf}) and transformations (\ref{CLN1}), we turn  $C^\prime\subset \mathcal{K}_\mathbf{z}[\Lambda^\prime]$ into a system of differential polynomials in $\mathcal{K}_\mathbf{z}[\Lambda,\xi_1^\prime]$ so that it connects the classical symmetries Z$(C^\prime)$ and nonclassical symmetries Z$(\mathcal{D})$.

From the structure of an ascending chain, we know that the both set $\mathcal{D}^{\prime\prime}$ and the complement set $C^{\prime}\backslash \mathcal{D}^{\prime\prime}$ are differential ascending chain in $\mathcal{K}_{\mathbf{z}}[\Lambda^\prime]$. Substituting (\ref{CLN1}) into $C^\prime$ and $\mathcal{D}^{\prime\prime}$, one obtains differential polynomial systems belonging to $\mathcal{K}_{z}[\Lambda, \xi_1^\prime]$ (still denote it as $C^\prime$ and $\mathcal{D}^{\prime\prime}$). It is noticed that a differential polynomial rank $\prec$ on $\mathcal{K}_\mathbf{z}[\Lambda^\prime]$ induces a corresponding rank on $\mathcal{K}_{\mathbf{z}}[\Lambda,\xi_1^\prime]$ through $\xi_i^\prime$ matching $\xi_i$ and $\eta^\prime$ matching $\eta (i=2, \cdots, p)$ and through keeping the highest rank of $\xi^\prime_1$. Under the induced rank, we do the reduction \begin{eqnarray}C=\textrm{Prem}((C^\prime\backslash \mathcal{D}^{\prime\prime})/\mathcal{D}^{\prime\prime})\subset \mathcal{K}_\mathbf{z}[\Lambda,\xi_1^\prime],\label{q}\end{eqnarray}
through formula (\ref{redf}). In the above reduction procedure for obtaining $C$, the reduction operation on $(\xi,\eta)$ is not involved. Therefore, due to the highest rank of $\xi_1^\prime$ and the linearity of both $C^\prime$ and transformations (\ref{CLN1}), each element in $C$, as a differential polynomial in $\mathcal{K}_\mathbf{z}[\Lambda, \xi_1^\prime]$ has the same rank and initials as its corresponding one in the $C^\prime\backslash \mathcal{D}^{\prime\prime}\subset \mathcal{K}_\mathbf{z}[\Lambda^\prime]$. Hence, the $C$ inherits the differential ascending chain property of the $C^\prime\backslash \mathcal{D}^{\prime\prime}$ as well as the linearity of its polynomials in their leading derivatives with (\ref{init}).

So far, we have finished the turning $C^\prime\subset \mathcal{K}_\mathbf{z}[\Lambda^\prime]$ into the $C\subset \mathcal{K}_\mathbf{z}[\Lambda,\xi_1^\prime]$.

For the sake of the following statement, we denote
the system of determining polynomials for a nonclassical symmetry as
\begin{equation*}
\mathcal{D}=\{p_1, p_2, \ldots, p_n\}\subset \mathcal{K}_\mathbf{z}[\Lambda],
\end{equation*}
and the set $C$ as
\begin{eqnarray*}
C=\{q_1, q_2, \ldots, q_m\} \subset \mathcal{K}_\mathbf{z}[\Lambda,\xi_1^\prime].\label{p1}
\end{eqnarray*}
For the obtained $C$, we have the connection with $C^\prime$ as follows.

\begin{theorem}\label{th1}
The set {\rm Z(}C{\rm )} contains all classical symmetries $\mathcal{X}^\prime=\Lambda^\prime \cdot \partial_\mathbf{z}$ with $\Lambda^\prime =(\xi_1^\prime, \cdots, \xi_p^\prime, \eta^\prime)$ and $\xi_1^\prime\not=0$ of  { \rm PDE (\ref{e1})} through transformations{ \rm (\ref{CLN1})}. More precisely, {\rm Z}$(C^\prime)\subset {\rm Z}(C)$ holds by means of the equivalence in definition 2.
\end{theorem}
Proof: For an $\Lambda^\prime\in \textrm{Z}(C^\prime)$, since a classical symmetry $\mathcal{X}^\prime=\Lambda^\prime \cdot \partial_\mathbf{z}$ with $\xi_1^\prime\not=0$ is equivalent to a nonclassical symmetry of the same PDE, there is a $\Lambda=(\xi_2, \cdots, \xi_p, \eta)\in {\rm Z}(D)$ connecting  $\mathcal{X}^\prime$ through transformations (\ref{CLN1}). The construction procedure for $C$ immediately shows that $\Lambda\in {\rm Z}(C)$. This ends the proof of the Theorem.

The relationship between $C$ and $\mathcal{D}$ is given in terms of the dependence of the determining polynomials in $\mathcal{D}$ on the ones in $C$ as stated in the following theorem.
\begin{theorem}\label{th3}
For each $p_i\in \mathcal{D}$, there exist differential operators $\mathfrak{D}_\nu^i$  with differential polynomial coefficients in $\mathcal{K}_\mathbf{z}[\Lambda]$ such that the identities
\begin{eqnarray}
{\rm IS}_i\cdot p_i=\sum_{\nu=1}^m \mathfrak{D}_\nu^i q_\nu, i=1, 2, \ldots, n,\label{key}
\end{eqnarray}
hold for all differential functions $\xi$, $\eta$ and an IS product ${\rm IS}_i$ of $C$.
\end{theorem}

Proof: From the construction of $C$, for each $q^\prime\in C^\prime\backslash D^{\prime\prime}$, there exist some differential operators $\mathfrak{D}_\nu^{\prime}$ and an IS product IS$^\prime$ of $D^{\prime\prime}$ such that
\begin{eqnarray*}
\textrm{IS}^\prime\cdot q^\prime=\sum_{d_\nu\in \mathcal{D}^{\prime\prime}}\mathfrak{D}_\nu^{\prime} d_\nu+q,
\end{eqnarray*}
from reduction formula (\ref{redf}). Here the remainder $q=\textrm{Prem}(q^\prime/\mathcal{D}^{\prime\prime})\in C$ is reduced with respect to $D^{\prime\prime}$. The obtained identity with assumption (\ref{init}) imply
\begin{eqnarray*}
\textrm{Z}(C)\cap \textrm{Z}(\mathcal{D}^{\prime\prime})\subset\textrm{Z}(C^\prime)\subset \textrm{Z}(C),
\end{eqnarray*}
from which we can easily prove that
\begin{eqnarray}
\textrm{Z}(C \cup \mathcal{D}^{\prime\prime})=\textrm{Z}(C^\prime). \label{zero1}
\end{eqnarray}
Furthermore, this equality shows that the differential ascending chain  $C\cup D^{\prime\prime}$ is irreducible since the $C^\prime$ is irreducible.

On the other hand, for each $p\in \mathcal{D}$, there are some differential operators $\mathfrak{D}_\nu$ and an IS product IS of $C$, such that
\begin{eqnarray}
\textrm{IS}\cdot p=\sum_{q_\nu\in C}\mathfrak{D}_\nu q_\nu+r, \label{pr1}
\end{eqnarray}
from reduction formula (\ref{redf}). Here the remainder $r=\textrm{Prem}(p/C)$ is reduced with respect to $C$.

Thus, identity (\ref{pr1}) and Theorem \ref{th1} imply that the remainder $r=\textrm{Prem}(p/C)$ is identically equal to zero on Z$(C^\prime)$, i.e.,
$$
r|_{\textrm{Z}(C^\prime)}=0.
$$
Furthermore, the equality (\ref{zero1}) indicates that
$$
r|_{\textrm{Z}(C\cup \mathcal{D}^{\prime\prime})}=0,
$$
which yields $r\equiv 0$ in (\ref{pr1}) due to $r$ being reduced with respect to $C\cup \mathcal{D}^{\prime\prime}$ and the irreducibility of $C\cup \mathcal{D}^{\prime\prime}$.
Consequently, the identities (\ref{key}) in the theorem are true. Hence this theorem is proved.

Thus, the identities (\ref{key}) with (\ref{init}) imply ${\rm Z}(C)\subset {\rm Z}(\mathcal{D}).$

Summarizing the above discussions, we have the inclusions:
\begin{eqnarray}\textrm{Z}(C^\prime)\subset \textrm{Z}(C)\subset \textrm{Z}(\mathcal{D}).\label{key2}
\end{eqnarray}

\textbf{Remark 5: }Theorems \ref{th2}-\ref{th3} and  Corollary 1 give  relationships between the classical and nonclassical symmetries of PDE (\ref{e1}) in terms of  the determining polynomial systems $C^\prime, C$ and $\mathcal{D}$ of the two kinds of symmetries, in which the set $C$ plays the rule of a bridge between classical and nonclassical symmetries of a PDE.

From the viewpoint of practical applications, the inclusions (\ref{key2}) imply that more additional equations from classical symmetry can be used to extended system $\mathcal{D}$ for nonclassical symmetry. In fact, any subset $Q \subset C$ is compatible with $\mathcal{D}$. Hence, we get some extended system of $\mathcal{D}$ by appending $Q$ to $\mathcal{D}$. The extended system $\mathcal{D}\cup Q$ is more easily solved than the original $\mathcal{D}$ since it contains more equations. It may yield a nonclassical symmetry of the underling PDE.

To summarize the above process of proving (\ref{key}), we obtain an algorithm for constructing (\ref{key}).

\textbf{Algorithm: }build the identities (\ref{key}).

\textbf{Input:} PDE (\ref{e1})

\textbf{Output:} Identities (\ref{key})

\textbf{Begin}

\hspace{0.3cm}\textbf{Step1. }Calculate the systems of determining equations $\mathcal{D}^\prime=0$ and $\mathcal{D}=0$ of the classical and nonclassical symmetries for a given PDE.

\hspace{0.3cm}\textbf{Step2. }Calculate the ascending chain form $C^\prime$ of $\mathcal{D}^\prime$ (by Wu's algorithm) and get $\mathcal{D}^{\prime\prime}$.

\hspace{0.3cm}\textbf{Step3. }Substitute (\ref{CLN1}) into $C^\prime$ and $\mathcal{D}^{\prime\prime}$ and calculate $C$ as in (\ref{q}).

\hspace{0.3cm}\textbf{Step4. }Establish (\ref{key}) by using formula (\ref{redf}).

\textbf{End.}

The following is an illustrative example to show how well the above theorems and algorithm work well.

\textbf{Example 1. }We consider the Burgers - Huxley equation
\begin{eqnarray}u_t=u_{xx}+u(u-1)(u-\sigma).\label{exam1}\end{eqnarray}

The whole classical symmetries admitted by the equation are obvious symmetries
\begin{eqnarray}
\mathcal{X}_1=\partial_t, \,\,\, \mathcal{X}_2=\partial_x,\label{CLBH}
\end{eqnarray}
corresponding to time and space translations for any $\sigma\in \mathds{R}$.

We first establish identities (\ref{key}) from the given algorithm.

\textbf{Step 1. }The systems of determining equations for a nonclassical symmetry $\mathcal{X}=\partial_t+\xi\partial_x+\eta\partial_u$ and a classical symmetry $\mathcal{X}^\prime=\tau^\prime \partial_t+\xi^\prime\partial_x+\eta^\prime\partial_u$ of the equation are given by $\mathcal{D}=\{p_1, p_2, p_3, p_4 \}=0$ and $\mathcal{D}^\prime=0$, where
\begin{eqnarray}
&&p_1=\xi_{uu},\nonumber\\
&&p_2=\eta_{uu}-2\xi_{xu}+2\xi\xi_u,\nonumber\\
&&p_3=2\eta_{xu}-\xi_{xx}+\xi_t+2\xi\xi_x-(3u(u-1)(u-\sigma)+2\eta)\xi_u,\label{Ncs}\\
&&p_4=\eta_{xx}-\eta_t+2\eta\xi_x+(\sigma-2u-2\sigma u+3u^2)\eta+u(u-\sigma)(u-1)(2\xi_x-\eta_u),\nonumber
\end{eqnarray}
and
\begin{eqnarray*}
&&\mathcal{D}^\prime=\left\{\tau^\prime_u, \tau^\prime_x, \tau^\prime_t -2 \xi^\prime_x, \xi^\prime_u, \xi^\prime_{xx}-2 \eta^\prime_{xu} -\xi^\prime_t, \eta^\prime_{uu},\right.\\
&&\left.\eta^\prime_t-\eta^\prime_{xx}+u(u-1)(u-\sigma ) (2\xi^\prime_x-\eta^\prime_u)+(\sigma +3 u^2-2 (\sigma+1) u)\eta^\prime\right\}.
\end{eqnarray*}

\textbf{Step 2. }Obviously, the system $\mathcal{\mathcal{D}}^\prime$ is already an ascending chain form under rank $t\prec x$ $\prec u\prec \xi^\prime\prec\eta^\prime\prec\tau^\prime$. Hence $C^\prime=\mathcal{D}^\prime$ and
\begin{eqnarray*}
\mathcal{D}^{\prime\prime}=\left\{\tau^\prime_t -2 \xi^\prime_x, \tau^\prime_x, \tau^\prime_u\right\}.
\end{eqnarray*}
with ${\rm{IS}}(C^\prime)=2\not \equiv 0$.

\textbf{Step 3.} Substituting the transformations (\ref{CLN1}) with $\xi_1^\prime=\tau^\prime$ into $C^\prime$, then computing the reduction Prem($(C^\prime\backslash\mathcal{D}^{\prime\prime})/\mathcal{D}^{\prime\prime}$) and after deleting the common factor $\tau^\prime\not=0$, we have the differential ascending chain $C=\{q_1, q_2, q_3, q_4\}$ under the induced rank $t\prec x\prec u\prec\xi\prec \eta$ of the rank in Step 2. Here
\begin{eqnarray}
&&q_1=\xi_{u},\nonumber\\
&&q_2=\eta_{uu},\nonumber\\
&&q_3=2\eta_{xu}+\xi_t+2\xi\xi_x-\xi_{xx},\label{Ccs}\\
&&q_4=\eta_{xx}-\eta_t+2\eta\xi_x+(\sigma-2u-2\sigma u+3u^2)\eta+u(u-\sigma)(u-1)(2\xi_x-\eta_u)\nonumber,
\end{eqnarray}
with IS$(C)=2.$

\textbf{Step 4. }Using reduction formula (\ref{redf}), one obtains relations (\ref{key}) as follows
\begin{eqnarray}
&&p_1=D_{u}q_1,\nonumber\\
&&p_2=q_2+2(\xi-D_x) q_1,\nonumber\\
&&p_3=q_3-(3u(u-1)(u-\sigma)+2\eta)q_1,\nonumber\\
&&p_4= q_4.\label{exam1c}
\end{eqnarray}

This implies Z$(C)\subset$Z$(D)$, i.e., the right side inclusion of (\ref{key2}) holds.

Obviously, the first classical symmetry in (\ref{CLBH}) with $\xi_1^\prime=\tau\equiv 1$ is contained in Z$(C)$, i.e., Theorem \ref{th1} is satisfied. We notice that the $(\xi,\eta)$ with
\begin{eqnarray*}
&&\xi(t,x,u)=a(x)=3 \tan(\frac{1}{2\sqrt{2}}(x+24 c_1))/{2 \sqrt{2}},\\
&&\eta(t,x,u)=(1-2 u)(8 a(x)^2+9)/48.
\end{eqnarray*}
is in Z$(C)$ for $\sigma=1/2$ but does not correspond to any classical symmetry. This shows that Z$(C)$ can be larger than the set of classical symmetries. The reason of this is that Z$(\mathcal{D}^{\prime\prime})=\emptyset$ for the current $\Lambda=(\xi,\eta)$ and equality (\ref{zero1}) is not true.

Rather, the $(\xi,\eta)$ with $\xi=(3u-\sigma-1)/\sqrt{2}, \eta=-(3/2)u(u-1)(u-\sigma)$ which makes all $p_i=0(i=1, 2, 3, 4)$ in (\ref{Ncs}), i.e., the $(\xi,\eta)\in {\rm Z}(\mathcal{D})$. However, it is obvious that $q_1\not=0$. Hence the $(\xi, \eta)\not\in{\rm Z}(C)$. This shows that Z$(C)$ is a proper subset of Z$(\mathcal{D})$ and the operator $\mathcal{X}=\partial_t+(3u-\sigma-1)/\sqrt{2}\partial_x-(3/2)u(u-1)(u-\sigma)\partial_u$ is a generator of a nontrivial nonclassical symmetry of equation (\ref{exam1}). Consequently, we properly have Z$(C^\prime)\subset {\rm Z}(C)\subset {\rm Z}(\mathcal{D})$.

\section{Applications}
We give two examples to show the applications of Theorems \ref{th2}-\ref{th3} and identities (\ref{key}) to determine nontrivial nonclassical symmetries of a given PDE.

\subsection{A nonclassical symmetry classification of a class of generalized Burgers equations}
The nonclassical symmetry classification of a PDE with arbitrary functions is a hard problem in symmetry analysis of a PDE \cite{D, I}. We solve this problem by using given Theorems \ref{th2}-\ref{th3} and identities (\ref{key}).

Consider the nonclassical symmetry classification of a class of generalized Burgers equations
\begin{eqnarray}
u_t+g(u)u_x-u_{xx}=f(u),\label{BH}
\end{eqnarray}
with two arbitrary functions parameters $f(u)$ and $g(u)$, which  arise in a wide range of mathematical models describing various processes in physics, biology and ecology(see \cite{C, RC1, RC2, RC3} and references therein).

For the case $g(u)=0$, the equation represents: the Huxley equation if $f(u)= u^2(1-u)$; the Fisher's equation if $f(u)=u(1-u)$; the generalized Fisher's equations if $f(u)=u-u^k$ or $f(u)=u^p(1-u^{p-1})$; the Chaffee -- Infante or Newell -- Whitehead equation if $f(u)=\beta u(1-u^2)$; the generalized KPP equation if $f(u)=\alpha u^3+\beta u^2+\rho u$; the Fitzhugh -- Nagumo equation if $f(u)=u(1-u)(u-\alpha)$, etc;

For the case $g(u)=\alpha u$, the equation represents: the Burgers equation if $f(u)=0$; the Burgers - Huxley equation if $f(u)=\beta u(1-u)(u-\gamma)$, etc;

For the case $g(u)=\alpha u^c (c\in R)$, the equation represents: the generalized Burgers -- Fisher equation if $f(u)=\beta u(1-u^c)$; the generalized Burgers -- Huxley equation if $f(u)=\beta u(1-u^c)(u^c-\gamma)$, etc.

In \cite{RC1}, the authors considered (\ref{BH}) with $g(u)=\lambda u$ and gave the nonclassical (there Q-conditional) symmetry classification with respect to $f(u)$. In \cite{RC2, RC3}, the complete classical symmetry classification of the equations was given.

Now we consider the nonclassical symmetry classification of (\ref{BH}) in general sense for two arbitrary parameters $f(u)$ and $g(u)$ by using the method given in the article.

The system of determining polynomials of classical symmetry $\mathcal{X}^\prime=\tau^\prime\partial_t+\xi^\prime\partial_x+\eta^\prime\partial_u$ of equation (\ref{BH}) is given by

$D^\prime=\left\{
\begin{array}{ll}
&\tau_u^\prime, \tau_x^\prime, \tau_t^\prime- 2 \xi^\prime_x, \xi^\prime_u,\eta^\prime_{uu},\\
&\eta^\prime_{xx}+f'(u)\eta^\prime+f(u)(2\xi^\prime_x-\eta^\prime_u)-g(u) \eta^\prime_x-\eta^\prime_t,\\
&2 \eta^\prime_{xu}+\xi^\prime_{xx}+g'(u)\eta^\prime+g(u)\xi^\prime_x-\xi^\prime_t.
\end{array}\right.$
\\
Under the order $t\prec x\prec u\prec \xi^\prime\prec\eta^\prime\prec\tau^\prime$, the system is already in differential chain form. Hence $C^\prime =\mathcal{D}^\prime$.

By using the given algorithm, we find the relationships (\ref{key}) between $\mathcal{D}$ and $C=\{q_1, q_2, q_3, q_4\}$  as follows
\begin{eqnarray}
&&p_1=q_3+2(\xi-g(u))q_4-2D_xq_4,\nonumber\\
&&p_2=q_2+(3f(u)-2\eta)q_4,\nonumber\\
&&p_3=D_{u}q_4,\nonumber\\
&&p_4=q_1,\label{classify}
\end{eqnarray}
where
\begin{eqnarray*}
&&q_1=f(u) \eta_u+g(u)\eta_x-\eta_{xx}+\eta_t+2( \eta-f(u)) \xi_x-\eta f^\prime(u),\\ &&q_2=2 \eta_{xu}-\xi_{xx}+(2 \xi-g(u)) \xi_x-\eta g^\prime(u)+\xi_t,\\
&&q_3=\eta_{uu},\\
&&q_4=\xi_u,
\end{eqnarray*}
with IS$(C)=-2$.

From the transformations (\ref{CLN1}) in this example, we have $C^\prime\ni\xi^\prime_u=\xi_u \tau^\prime$ which implies that equations (\ref{BH}) admits a nonclassical symmetry when $q_4\not=0$ and Z$(\mathcal{D})\not=\emptyset$ from Theorem \ref{th2}. Therefor, based on the information we consider the cases of $q_4\not=0$  and  $q_4=0$ respectively.

{\bf Case 1: $q_4\not=0.$}

In the case, from Theorem \ref{th2} we can impose more restrictions on the system $\mathcal{D}=0$ so that it is solved easily and yields nontrivial nonclassical symmetries of the equation. It is  observed that if one takes $\eta$ as $\eta_{uu}\not=0$, then $\eta_{uu}^\prime=\eta_{uu}\tau^\prime\in C^\prime$ never be zero under transformation (\ref{CLN1}) for $\tau^\prime\not=0$. Hence, according to Theorem \ref{th2}, such $\xi$ and $\eta$ may yield a nontrivial nonclassical symmetry of the equation. Based on the analysis, an ansatz of the infinitesimal functions that fits these conditions is $\xi=a(t,x)u+b(t,x)$ and $\eta=\alpha(t,x)u^3+\beta(t,x)u^2+\rho(t,x)+\gamma(t,x)$ with $a^2(t,x)+\beta^2(t,x)\not=0$ or $\alpha(t,x)\not=0$.

\textbf{Remark 6:} Actually, one takes $\eta=\sum_{i=0}^n\alpha_i(t,x)u^i$ in more general form. However, it is easy to deduce that $\alpha_i\equiv 0$ when $i\geq 4$.

Under these restrictions, system $\mathcal{D}=0$ is reduced significantly and solved easily. In particular, we have the following three groups of  nonclassical symmetries.
\begin{eqnarray}\xi=a u+b, \eta=\frac{1}{3}(c_1-a)u^3+a(c_2-b)u^2+\rho u+\gamma,\label{s1}\end{eqnarray}
for $f(u)=\frac{\left(2 a+c_1\right)}{9 a} \left(a (c_1-a)u^3+3 a (c_2-b)u^2+3 \rho  u+3 \gamma \right), g(u)=c_1 u+c_2$ where $ a, b, c_1, c_2, \rho$ and $\gamma$ are arbitrary constants with $a\not=0$ or $c_2\not=b$ or $c_1\not=a$;
\begin{eqnarray}\xi=\frac{c_1}{4}u+b, \eta=\frac{1}{16} c_1^2 u^3-\frac{1}{4} c_1(b-c_2) u^2+\rho  u+\gamma,\label{s2}\end{eqnarray}
for $f(u)=\frac{1}{8} c_1^2 u^3+\frac{1}{2} c_1(c_2-b) u^2+2 \rho  u+2 \gamma, g(u)=c_1u+c_2$, where $c_1, c_2, b, \rho$ and $\gamma$ are arbitrary constants with $c_1\not =0$;
\begin{eqnarray}\xi=-\frac{c_1}{2}u+b(t,x),\eta=-\frac{1}{4} c_1^2 u^3+\frac{1}{2} c_1 u^2 (b(t,x)-c_2)+ \rho (t,x)u+\gamma (t,x),\label{s3}\end{eqnarray}
for $f(u)=0, g(u)=c_1u+c_2$, where $b=b(t,x), \rho=\rho(t,x)$ and $\gamma=\gamma(t,x)$ are solutions to the system
\begin{eqnarray*}
&&b_t-b_{xx}+(2b-c_2)b_x-2 \rho_x=0,\\
&&\rho_t-\rho_{xx}+c_2 \rho_x+ 2 \rho b_x+c_1 \gamma_x=0,\\
&&\gamma_t -\gamma_{xx}+c_2 \gamma_x+2 \gamma b_x=0,
\end{eqnarray*}
where $c_1$ and $c_2$ are arbitrary  constants with $c_1\not=0$.

\textbf{Case 2: $q_4=0$.}

In the case, the condition $q_4=\xi_u=0$ gives a first step simplification of the system $\mathcal{D}$ of determining equations. In order to further simplify the system, we use Theorem \ref{th2} to find more information on nontrivial nonclassical symmetry of the equation. For a $\Lambda=(\xi,\eta)\in $Z$(\mathcal{\mathcal{D}})$, under transformations (\ref{CLN1}), we have
$$D^{\prime\prime}=\{\tau^\prime_u,\tau^\prime_x, \tau^\prime_t-2\xi_x\tau^\prime\}.$$

It is clear that Z$(\mathcal{D}^{\prime\prime})=\emptyset$ if $\xi_{xx}\not=0$. Hence, by  Corollary 1, we look for a nontrivial nonclassical symmetry of the form $\xi=a(x)$ with $a^{\prime\prime}(x)\not =0$. Now, we suppose $\eta=b(t,x)u+c(t,x)$ from the $p_1=0$. These significantly simplify the system of determining equations.

For simplicity, denote $b=b(t, x), c=c(t, x)$ and $a=a(x)$.

Now, the rests of determining equations $p_2=0$ and $p_4=0$ become
\begin{eqnarray}
&&(u+A)f^\prime(u)+(2B-1)f(u)+(E u+F)g(u)+H u=G,\nonumber\\
&&(u+A)g^\prime(u)+Bg(u)=K,\label{clas0}
\end{eqnarray}
for $b\not=0$, where
\begin{eqnarray*}
&&A=c/b, B=a^\prime/b, K=(2aa^\prime+2b_x-a^{\prime\prime})/b, H=(b_{xx}-b_t-2ba^\prime)/b,\\
&&E=-b_x/b, F=-c_x/b, G=(c_t-c_{xx}+2ca^\prime)/b;
\end{eqnarray*} and
\begin{eqnarray}
&&f^\prime(u)+2B_1f(u)+K_1g(u)=H_1,\nonumber\\
&&g^\prime+B_1g(u)=E_1\label{clas1}
\end{eqnarray}
for $b=0$ and $c\not=0$, where
\begin{eqnarray*}B_1=a^\prime/c, K_1=-c_x/c, E_1=(2aa^\prime-a^{\prime\prime})/c, H_1=(2ca^\prime +c_t-c_{xx})/c.\end{eqnarray*}

Since the $f$ and $g$ are functions of a single variable $u$, the coefficients $A, B, K,$ $ H, E, F, G, , B_1, E_1, K_1$ and $H_1$  must be constants. After simple analysis, we can prove that equation (\ref{BH}) does not admit a nontrivial nonclassical symmetry unless $g(u)\equiv 0, b=b(x)\not=0$  and $c=c(x)$.

Under $g(u)=0$, equations (\ref{clas0}) become
\begin{eqnarray}
&&(u+A)f^\prime(u)+(2B-1)f(u)+H u=G,\nonumber\\
&&K=2aa^\prime+2b^\prime-a^{\prime\prime}=0,\label{clas01}
\end{eqnarray}
with
\begin{eqnarray}c=Ab,a^\prime=Bb, b^{\prime\prime}-2a^\prime b=Hb, 2a^\prime c-c^{\prime\prime}=Gb, AH+G=0.\label{coef}\end{eqnarray}

Solving the first equation, we get
\begin{eqnarray}f(u)=k (u+A)^\mu-\frac{H}{2 B}(u+A),\textrm{ with }\mu=1-2B\not=0,\label{f1}\end{eqnarray}
and
\begin{eqnarray}f(u)=k-H u,\textrm{ with }2B-1=0,\label{f2}\end{eqnarray}
for arbitrary constant $k$.

It easily proves that the compatibility condition of (\ref{coef}) for case (\ref{f1}) of $f(u)$ is $B=-1$, i.e., $\mu=3$. Meantime, the $a=a(x)$ satisfies $3a^{\prime\prime}-2aa^\prime=0$ which further determines $b$ and $c$ from theequations in (\ref{coef}). Consequently, we obtain nonclassical symmetries of equation (\ref{BH}) as follows.
\begin{eqnarray}
&&\xi=a(x),\,\,\, \eta=(u+A)b(x),\label{CLK0}
\end{eqnarray}
with $b(x)=-a^\prime(x), c(x)=-A a^\prime(x)$ and
\begin{eqnarray}
&&a(x)=k_1\tanh((x+3c_1)k_1/3),\textrm{ for } f(u)=k(u+A)^3+H(u+A)/2;\nonumber\\
&&a(x)=-3/(x+3c_2),\textrm{ for }f(u)=k(u+A)^3,\label{Clk0}
\end{eqnarray}
where $c_1$ and $c_2$ are arbitrary constants and $k_1^2=-9H/4\not=0$.

In case (\ref{f2}), the equation (\ref{BH}) becomes a linear equation.
The determining equations for nonclassical symmetries of the linear equation yield
\begin{eqnarray}
&&\xi=c(t,x),\,\,\, \eta=a(t,x)u+b(t,x),\label{CLK1}
\end{eqnarray}
where the functions $a=a(t,x), b=b(t,x)$ and $c=c(t,x)$ satisfy
\begin{eqnarray*}
&&b_t-b_{xx}+2(b-k)c_x+k a+Hb=0,\\
&&a_t-a_{xx}+2 (a+H)c_x=0,\\
&&c_t-c_{xx}+2 cc_x+ 2 a_x=0.
\end{eqnarray*}
It is not possible to find the general solutions of the above PDE system. In the following, we give some special solutions.

Solving the system under the condition $b(t,x)=(-k/H)a(t,x)$, we obtain two groups of solutions:
\begin{eqnarray}
&&a(t,x)=\frac{2 c_2^2(c_1+2 c_2 c_4) \tanh (T)+2c_2(c_4 c_1-2c_2 H)+c_1^2}{4c_2^2}, \nonumber\\
&&c(t,x)=c_4-c_2 \tanh (T);\label{news1}
\end{eqnarray}
and
\begin{eqnarray}
&&a(t,x)=\frac{4c_2^4-4Hc_2^2-c_1^2-6c_2^2(2c_2^2\tanh(T)^2+c_1\tanh(T))}{4c_2^2},\nonumber\\
&&c(t,x)=-\frac{3 c_2^2 \tanh (T)+c_1}{c_2},\label{news2}
\end{eqnarray}
with $T=c_1 t+c_2 x+c_3$ in terms of arbitrary constants $c_1,c_2\not=0$ and $c_3$.

In the special case where the functions $a=a(x), b=b(x)$ and $c=c(x)$ are independent of $t$, we have a additional group of solutions to the determining equations given by
\begin{eqnarray}
&&a(x)=\frac{c^\prime(x)-c^2(x)}{2}-H, b(x)=-\frac{k}{H} a(x), c(x)=\frac{3x^2}{3k_1-x^3},\label{Clk1}
\end{eqnarray}
in terms of an arbitrary constants $k_1$.

Summarizing the above procedure, we obtain a nonclassical symmetry classification of equation (\ref{BH}). We have proved that the equation  admits nontrivial nonclassical symmetry when $f(u)$ is a cubic (including linear) and $g$ is linear in $u$. The symmetries (\ref{CLK0}) with (\ref{Clk0}) and (\ref{CLK1}) with (\ref{Clk1})  were given in \cite{C}. The symmetries (\ref{s1}), (\ref{s2}), (\ref{CLK1}) with (\ref{news1}) and (\ref{news2}) are new ones. The symmetries (\ref{s3}) are more general than those given in \cite{C}.

Moreover, here in the case $g=0$,  we recovered the symmetry classification results in \cite{C} obtained by P. A. Clarkson and E. L. Mansfield who directly treat the system of determining equations by the Gr${\rm\ddot{o}}$ebner basis method. In the case $g(u)=\lambda u$, we recovered the results in \cite{RC1} given by Roman Cherniha who dealt with the determining system by direct calculation. Here, we solve the problem for equation (\ref{BH}) by an alternative method using Theorems \ref{th2}-\ref{th1} and identities (\ref{key}) in general sense for arbitrary $f(u)$ and $g(u)$ and show that only in the case when $g$ is linear and $f(u)$ is cubic polynomials in $u$, the class of equations (\ref{BH}) admits a nontrivial nonclassical symmetry.  In particular, compared with the ways for symmetry classification problem existing in literatures\cite{C, D, chao5, RC1}, our method avoids a lot calculations for dealing with the determining equations.

Of physical interest is that there are many important equations  such as the Huxley equation, the Chaffee - Infante or Newell - Whitehead equation, the generalized Fisher's equations with $k=3$ or $p=2$, the generalized KPP equation, the Fitzhugh - Nagumo equation, Burgers equation, the Burgers - Huxley equation, etc., admit nontrivial nonclassical symmetries. Nevertheless, equations such as the generalized Burgers-Fisher equation, the generalized Burgers - Huxley equation ($c\not=1$), the Fisher's equation, the generalized Fisher's equations with $k\not=3$ and $p\not=2$, etc., do not admit nontrivial nonclassical symmetries.

\subsection{Nonclassical symmetry of a KdV-type equations}
Consider the coupled system of KdV-type systen given by
\begin{eqnarray} u_t+3vv_x=0,\,\,\,
v_t+2v_{xxx}+2uv_x+u_xv=0,\label{e23}
\end{eqnarray}
which are found in the study of Kac-Moody algebras and in soliton theory \cite{Gerd}.

The system of determining polynomials for a classical symmetry of the equations is given by
\begin{eqnarray*}\mathcal{D}^\prime=\left\{
\begin{array}{ll}
&\tau_v^\prime, \tau_u^\prime, \tau_x^\prime, 3 \eta^\prime+2 u \tau_t^\prime,\\
&\eta_v^\prime, \eta_x^\prime, \eta_t^\prime, \eta^\prime-u \eta_u^\prime, u \phi^\prime-v \eta^\prime,\\
&\xi_v^\prime, \xi_u^\prime, \xi_t^\prime, \eta^\prime+2 u \xi_x^\prime.
\end{array}\right.
\end{eqnarray*}
and it has ${\rm Z}(\mathcal{D}^\prime)=\{\tau^\prime, \xi^\prime, \eta^\prime, \phi^\prime\}$, where
\begin{eqnarray*}
\tau^\prime=c_2-3c_1 t,\,\, \xi^\prime=c_3-c_1 x, \,\, \eta^\prime=2c_1 u,\,\, \phi^\prime=2c_1 v,
\end{eqnarray*}
and $ c_1, c_2, c_3\in \mathds{R}$ are arbitrary constants. It yields obvious classical symmetries
\begin{eqnarray*}
\mathcal{X}^\prime_1=-\frac{1}{2}(3t\partial_t+x\partial_x)+u\partial_u+v\partial_v, X^\prime_2=\partial_t,\mathcal{X}^\prime_3=\partial_x
\end{eqnarray*}
of equations (\ref{e23}). By our given algorithm, we have the set $C=\{q_1, \cdots, q_9\}$, where
\begin{eqnarray}
\begin{array}{ll}
&q_1=3\eta\xi-2u\xi_t,\\
&q_2=2u\xi_x+\eta,\\
&q_3=u\phi-v\eta,
\end{array}
\begin{array}{ll}
&q_4=3\eta^2 -2u\eta_t,\\
&q_5=\eta-u\eta_u,\\
&q_6=\xi_u, q_7=\xi_v, q_8=\eta_x, q_9=\eta_v.
\end{array}
\end{eqnarray}

The system $\mathcal{D}$ of determining polynomials of the nonclassical symmetry of the equations is consists of 11 strongly nonlinear equations, variable coefficients and long expressions for differential polynomials. We can establish the identities (\ref{key}) between the system $\mathcal{D}$ and $C$ by our given algorithm. Due to complexity and massive scale, we here omit the expressions of the connections. The system $\mathcal{D}$ is difficult to deal with directly. Nevertheless, we can get some specific nonclassical symmetries of equation (\ref{e23}) by using identities (\ref{key}) and inclusions (\ref{key2}).

We notice that there are three identities in (\ref{key}) that only involve the $q_6$ and $q_7$. In order to reduce the system, we deliberately set $q_6=q_7=0$. In addition, we look for a projection symmetry by setting $q_9=0$ ($\eta$ does not involve coordinate $v$) to further simplifying the system $\mathcal{D}=0$. Let $Q=\{q_6, q_7, q_9\}\subset C$ (see the statement below Remark 5). For the extended polynomial system $\mathcal{D}\cup Q$, using Wu's algorithm in \cite{wu3,chao5},  we obtain  a differential ascending chain (actually it is a characteristic set)
\begin{eqnarray}
\mathcal{AS}=\left\{\begin{array}{l} \eta_v,\phi_u,\xi_v,\xi_u,\phi_x,\phi-v \phi_v, 2 v (\xi-u) \xi_x-v \eta +\xi \phi,\\
v(\xi-u)\eta_u+v \eta-\xi\phi, 2 v(\xi-u)\eta_t-v \eta^2-(3 \xi-4 u)\phi\eta,\\
2v(\xi-u)^2\phi_t-(3\xi-4u)\xi\phi^2+2v(\xi-2u)\phi+(v\eta)^2,\\
2v^2 (\xi-u)^2\eta_x -v^2\eta^2+v(\xi +u)\eta\phi-u \xi\phi^2,\\
2 v (\xi-u)\xi_t-((3 \xi- 4 u) \phi+v\eta) \xi, \end{array} \right\}\nonumber
\end{eqnarray}
of extended system $\mathcal{D}\cup Q$.

Solving the system, we find that Z$(\mathcal{AS})=\{\xi,\eta,\phi\}$ where
\begin{eqnarray*}
&&\begin{array}{ll}
\xi=\frac{1}{3} x \left(\tanh \left(T\right)+1\right) c_1+ e^{c_1 t}  \text{sech}\left(T\right)c_3,\\
\eta= \left(\frac{1}{3} x \left(\tanh \left(T\right)+1\right)c_1+ e^{c_1 t} \text{sech}\left(T\right)c_3-\frac{2}{3} u \left(\tanh
   \left(T\right)+1\right)\right)c_1,\\
\phi=\frac{1}{3} v \left(1-2 \tanh \left(T\right)\right)c_1;\\
\end{array}\\
&&
\begin{array}{ll}
\xi=\frac{1}{3} x \left(\coth \left(T\right)+1\right)c_1+e^{c_1 t} \text{csch}\left(T\right)c_3 ,\\
\eta=\left(\frac{1}{3} x \left(\coth \left(T\right)+1\right)c_1+e^{c_1 t} \text{csch}\left(T\right)c_3 -\frac{2}{3}u \left(\coth
   \left(T\right)+1\right)\right)c_1,\\
\phi=\frac{1}{3} v \left(1-2 \coth \left(T\right)\right)c_1;\\
\end{array}\\
&&
\begin{array}{ll}
\xi=c_2 e^{2c_1t},\\
\eta=c_1 c_2 e^{2c_1t},\\
\phi=c_1 v,\\
\end{array}
\end{eqnarray*} in which $T=c_1 t-c_2$ and $c_1(\not=0), c_2, c_3$ are arbitrary constants. Obviously, these solutions are not in Z$(\mathcal{D}^\prime)$. Hence the coupled KdV-type system (\ref{e23}) admits nontrivial nonclassical symmetries with infinitesimal functions given above.

\section{Conclusion remarks}
Some relationships (Theorems \ref{th2}-\ref{th3}) between the classical  and nonclassical symmetries of a PDE are given. A interesting point is that by the equivalent transformation (\ref{CLN1}) infinitesimal generators of classical and nonclassical symmetries of a PDE, we derive a set $C$ which connect such symmetries of the PDE in terms of their determining polynomials. In addition, Ritt -- Wu's reduction formula and algorithm play fundamental rules in deriving these results. In particular, an algorithm based on the reduction algorithm is given for constructing the connection (\ref{key}). These connections efficiently yield extra conditions on the nonclassical symmetries of a PDEs, which make it easier to solve the system of nonlinear determining equations. Since the connections are established by using a pseudo-reduction formula for a differential polynomial with respect to the ascending chain of the system of determining polynomials of a classical symmetry, it is not necessary to compute a characteristic set of the system. Hence the computation is much more cheaper. Consequently, the connections give an alternative way to obtain nontrivial nonclassical symmetry of a PDE and significantly improve the efficiency of the algorithms given in literatures \cite{D, chao5}. Moreover, the connections yield a way to better understand the relationship between classical and nonclassical symmetries of given PDEs.  The higher efficiency applications of the obtained results for determining nontrivial nonclassical symmetries yielded new nonclassical symmetries for both a class of the generalized Burgers equations and a coupled KdV type system.

\textbf{Acknowledgments: }{The first author are supported by Natural science foundation of China, under  grand numbers: 11571008}. George Bluman acknowledges financial support from the National Sciences and Engineering Research Council of Canada.

\end{document}